\documentclass[11pt,letterpaper]{amsart}

\usepackage{amssymb} 
\usepackage{amsthm}
\usepackage{amsfonts}
\usepackage{amsmath}
\usepackage{amstext}
\usepackage{amscd}

\numberwithin{equation}{section}

\theoremstyle{plain}

\newtheorem{thm}{Theorem}[section]
\newtheorem{prop}[thm]{Proposition}
\newtheorem{cor}[thm]{Corollary}

\theoremstyle{definition}

\theoremstyle{remark}
\newtheorem{rem}[thm]{Remark}
\newtheorem{question}{Question}

\newcommand{\one}{\mathbf {1}}
\newcommand{\cons}{\operatorname{F }}

\newcommand{\sign}{\mbox{$\operatorname {sgn}$}\,}
\DeclareMathOperator{\hlink}{\tilde \Lambda }
\DeclareMathOperator{\link}{ \Lambda }

\DeclareMathOperator{\ahlink} {\tilde \Omega }

\newcommand{\C}{\mathbb{C}}
\newcommand{\R}{\mathbb{R}}

\newcommand{\Z}{\mathbb{Z}}
\newcommand{\N}{\mathbb{N}}
\newcommand{\Q}{\mathbb{Q}}
\newcommand{\e}{\varepsilon}

\newcommand{\inv}{^{-1}}

\newcommand{\cX}{\overline {X}}

\newcommand{\half}{\tfrac12} \def\({(\!(} \def\){)\!)}
 \newcommand{\dchi}{\, d\chi}

\date{February 7, 2002}

\thanks{Research supported by CNRS and NSF grant DMS-9972094}

\title[Algebraically Constructible Functions] {Algebraically
Constructible Functions:\\ Real Algebra and Topology}

\author{Clint McCrory and Adam Parusi\'nski}
\address{Department of Mathematics, University of Georgia,
Athens, GA 30602, USA }
\email{clint@math.uga.edu}
\address{D\' epartement de Math\'ematiques, Universit\'e
d'Angers, 2, bd Lavoisier, 49045 Angers Cedex 01, France}
\email{parus@tonton.univ-angers.fr}

\subjclass{Primary: 14P25. Secondary: 14B05, 14P10}

\newcommand\abstracttext{}

\begin{document} \begin{abstract} \abstracttext {Algebraically
constructible functions connect real algebra with the topology
of algebraic sets. In this survey we present some history,
definitions, properties, and algebraic characterizations of
algebraically constructible functions, and a description of
local obstructions for a topological space to be homeomorphic to
a real algebraic set.} \end{abstract} \maketitle


More than three decades ago Sullivan proved that the link of
every point in a real algebraic set has even Euler
characteristic. Related topological invariants of real algebraic
singularities have been defined by Akbulut and King using
resolution towers and by Coste and Kurdyka using the real
spectrum and stratifications.

Sullivan's discovery was motivated by a combinatorial formula
for Stiefel-Whitney classes. Deligne interpreted these classes
as natural transformations from constructible functions to
homology. Constructible functions have interesting operations
inherited from sheaf theory: sum, product, pullback,
pushforward, duality, and integral. Duality is closely related
to a topological \emph{link operator.} To study the topology of
algebraic sets the authors introduced \emph{algebraically
constructible functions}. Using the link operator we have
defined many local invariants which generalize those of
Akbulut-King and Coste-Kurdyka.

Algebraically constructible functions are interesting from a
purely algebraic viewpoint. From the theory of basic algebraic
sets it follows that if a constructible function $\varphi$ on an
algebraic set of dimension $d$ is divisible by $2^d$ then
$\varphi$ is algebraically constructible. Parusi\' nski and
Szafraniec showed that algebraically constructible functions are
precisely those constructible functions which are sums of signs
of polynomials. Bonnard has given a characterization of
algebraically constructible functions using fans, and she has
investigated the number of polynomials necessary to represent an
algebraically constructible function as a sum of signs of
polynomials. Pennaneac'h has developed a theory of \emph{algebraically
constructible chains} using the real spectrum.

In section 1 we briefly discuss the results of Sullivan,
Akbulut-King, and Coste-Kurdyda. In the next section we define
algebraically constructible functions and their operations. In
section 3 we discuss the relations of algebraically
constructible functions with real algebra. In the following
section we describe how to generate our local topological
invariants. In the final section we raise some questions for
future research. Throughout we consider only algebraic subsets
of affine space.

Related survey articles have been written recently by Coste
\cite{cos2}, Bonnard \cite{bon2}, and McCrory \cite{mcc}. We
thank Michel Coste for his encouragement and insight.


\medskip \section{Akbulut-King Numbers}\label{ak}

Let $X$ be a real semialgebraic set in $\R^n$, and let $x\in X$.
Let $S(x,\e)$ be the sphere of radius $\e > 0$ in $\R^n$
centered at $x$. By the local conic structure lemma \cite{bcr}
(9.3.6), for $\e$ sufficiently small the topological type of the
space $S(x,\e)\cap X$ is independent of $\e$. This space is
called the \emph{link} of $x$ in $X$, and it is denoted by
$\operatorname{lk}(x,X)$.

Our starting point is Sullivan's theorem \cite{sul}:

\begin{thm}\label{sullivan} If $X$ is a real algebraic set in
$\R^n$ and $x\in X$ then the Euler characteristic
$\chi(\operatorname{lk}(x,X))$ is even.\end{thm}

For example, the ``theta space'' $X\subset \R^2$,
\begin{equation*}X = \{(x,y)\ |\ x^2+y^2 = 1\}\cup \{(x,y)\ |
-1\leq x \leq 1,\ y = 0 \},\end{equation*} is not homeomorphic
to an algebraic set, for the link of the point $(1,0)$ (or the
point $( -1,0)$) in $X$ is three points, which has odd Euler
characteristic.

Many proofs of Sullivan's theorem have been published; see
\cite{burver}, \cite{har}, \cite{benris} (3.10.4), \cite{bcr}
(11.2.2), \cite{fumcc} (4.4). Sullivan's original idea was to
use complexification. First he proved that the link of $x$ in
the complexification $X_\C$ has Euler characteristic 0, and then
he used that $\operatorname{lk}(x,X)$ is the fixed point set of
complex conjugation on $\operatorname{lk}(x,X_\C)$ to deduce that
\begin{equation*}\chi(\operatorname{lk}(x,X))
\equiv\chi(\operatorname{lk}(x,X_\C))
\pmod 2.\end{equation*} Mather \cite{mat} (p.~221) gave a proof
that the link $L$ of a point in a complex algebraic set has
Euler characteristic 0 by constructing a tangent vector field on
$L$ which integrates to a nontrivial flow of $L$.

The following result puts Sullivan's theorem in a more general
context (\emph{cf.} \cite{akbking2} (2.3.2)).

\begin{thm}\label{family} If $X$ and $Y$ are real algebraic sets
with $Y$ irreducible and $f:X\to Y$ is a regular map, there is
an algebraic subset $Z$ of $Y$ with $\dim Z < \dim Y$ such that
the Euler characteristic $\chi(f\inv(y))$  is constant mod 2 for
$y\in Y\setminus Z$.\end{thm}

In other words, the Euler characteristic is generically constant
mod 2 in every family of real algebraic sets.  To deduce
Sullivan's theorem as a corollary let $Y = \R$, $x_0\in X$, and
$f(x) = (x-x_0)^2$. For $y<0$ the fiber $f\inv(y)$ is empty, and
for $y > 0$ sufficiently small, the fiber $f\inv(y)$ is
$\operatorname{lk}(x_0,X)$.

Benedetti-Ded\`o \cite{bended} and Akbulut-King \cite{akbking1}
proved that Sullivan's condition is not only necessary but also
sufficent in low dimensions: If $X$ is a compact triangulable
space of dimension less than or equal to 2, and the link of
every point has even Euler characteristic, then $X$ is
homeomorphic to a real algebraic set.  (The link of a point in a
triangulable space is the boundary of a simplicial
neighborhood.) A triangulable space such that the link of every
point has even Euler characteristic is called an \emph{Euler
space}.

Akbulut and King \cite{akbking2} showed that the situation in
dimension 3 is more complicated. They defined four non-trivial
topological invariants of a compact Euler space $Y$ of dimension
at most 2, $a_i(Y)\in \Z/2$, $i=0,1,2,3$ (with
$a_i(Y) = 0$ when $\dim Y < 2$). Let $\chi_2(Y)$ be the Euler
characteristic mod 2. It is easy to see that if $X$ is an Euler
space then the link of every point of $X$ is an Euler space.

\begin{thm}\label{akbulutking} A compact 3-dimensional
triangulable topological space $X$ is homeomorphic to a real
algebraic set if and only if, for all $x\in X$,
$\chi_2(\operatorname{lk}(x,X)) = 0$ and
$a_i(\operatorname{lk}(x,X)) = 0$, $i = 0,1,2,3.$ \end{thm}

Akbulut and King's invariants arise from a combinatorial
analysis of the resolution of singularities of an algebraic set.
The elementary definition of these \emph{Akbulut-King numbers}
and computations of examples can be found in Akbulut and King's
monograph \cite{akbking2}, chapter VII, pages 190--197. (In the
terminology of \cite{akbking2} (7.1.1),
$a_i(\operatorname{lk}(x,X))$ is the mod 2 Euler characteristic
of the link of $x$ in the characteristic subspace $\mathcal
Z_i(X)$.) The depth of the method of \emph{resolution towers} is
shown by the remarkable result that the vanishing of the
Akbulut-King numbers gives a sufficient condition for a
triangulable 3-dimensional space to be homeomorphic to an
algebraic set. Chapter I of \cite{akbking2} is an introduction
to their methods, with informative examples.

Another descendant of Sullivan's theorem is due to Coste and
Kurdyka \cite{coskur}:

\begin{thm}\label{costekurdyka} Let $X$ be an algebraic set and
let $V$ be an irreducible algebraic subset. For $x\in V$ the
Euler characteristic of the link of $x$ in $X$ is generically
constant mod 4: There is an algebraic subset $W$ of $V$ with
$\dim W < \dim V$ such that $\chi(\operatorname{lk}(x,X))$ is
constant mod 4 for $x\in V\setminus W$. \end{thm}

This theorem was first proved by Coste \cite{cos1} when $\dim X
- \dim V \leq 2$ using chains of specializations of points in
the real spectrum. The general case was proved topologically
using stratifying families of polynomials. It can also be proved
using Akbulut and King's topological resolution towers
(\cite{akbking2}, exercise on p.~192).

Using the same techniques Coste and Kurdyka defined invariants
mod $2^k$ associated to chains $X_1\subset X_2\subset \cdots
\subset X_k$ of algebraic subsets of $X$ (\cite{coskur}, Theorem
4). Furthermore they used their mod 4 and mod 8 invariants to
recover the Akbulut-King numbers. Using a relation between
complex conjugation and the monodromy of the complex Milnor
fibre of an ordered family of functions, the authors
\cite{mccpar1} reinterpreted and generalized the Coste-Kurdyka
invariants as Euler characteristics of iterated links.


\medskip \section{Constructible Functions}\label{cons}

Algebraically constructible functions were introduced by the
authors \cite{mccpar2} as a vehicle for using the ideas of Coste
and Kurdyka to generate new Akbulut-King numbers.

Let $X$ be a real semialgebraic set.  A \emph{constructible
function} on $X$ is an integer-valued function
\begin{equation*}\varphi:X\to \Z\end{equation*} which can be
written as a finite sum \begin{equation}\label{cf} \varphi =
\sum m_i \one_{X_i},\end{equation} where for each $i$, $X_i$ is
a semialgebraic subset of $X$, $\one_{X_i}$ is the
characteristic function of $X_i$, and $m_i$ is an integer.

The set of constructible functions on $X$ is a ring under
pointwise \emph{sum} and \emph{product}. If $f:X\to Y$ is a
semialgebraic map and $\varphi$ is a constructible function on
$Y$, the \emph{pullback} $f^*\varphi$ is the constructible
function defined by
\begin{equation}\label{pullback}f^*\varphi(x) =
\varphi(f(x)).\end{equation}

The operations of pushforward and duality are defined using the
Euler characteristic. If $\varphi$ has compact support one may
assume that the sets $X_i$ in (\ref{cf}) are compact, and the
\emph{Euler integral} is defined by
\begin{equation}\label{integral}\int_X\varphi\ \dchi = \sum
m_i\chi(X_i).\end{equation} The Euler integral is additive, and
it does not depend on the choice of representation (\ref{cf}) of
$\varphi$.

If $f:X\to Y$ is a proper semialgebraic map and $\varphi$ is a
constructible function on $X$, the \emph{pushforward},
$f_*\varphi$ is the constructible function on $Y$  given by
\begin{equation}\label{pushforward} f_*\varphi (y) =
\int_{f^{-1}(y)} \varphi\ \dchi.\end{equation}

Suppose that $X$ is a semialgebraic set in $\R^n$.  If
$\varphi$ is a constructible function on $X$, the \emph{link}
$\link\varphi$ is the constructible function on $X$ defined by
\begin{equation}\label{link}\Lambda\varphi (x) =
\int_{S(x,\varepsilon)\cap X} \varphi\ \dchi,\end{equation} for
$\varepsilon>0$ sufficiently small.

The \emph{dual} $D\varphi$ is defined by
\begin{equation}\label{dual}D\varphi = \varphi -
\link\varphi.\end{equation}

The operations sum, product, pullback, pushforward, and dual
come from sheaf theory. Operations on constructible functions
have been studied by Kashiwara and Schapira \cite{kassch}
\cite{sch} and by Viro \cite{vir}.

Now suppose that $X$ is a real algebraic set. A provisional
definition of algebraically constructible functions would be to
require the sets $X_i$ in (\ref{cf}) to be algebraic subsets of
$X$. But the image of an algebraic set by a proper regular map
is not necessarily algebraic, so this class of functions---which
we call \emph{strongly algebraically constructible}---is not
preserved by the pushforward operation. To remedy this defect we
make the following definition.

Let $X$ be a real algebraic set. An \emph{algebraically
constructible function} on $X$ is an integer-valued function
which can be written as a finite sum
\begin{equation}\label{acf}\varphi = \sum m_i
f_{i*}\one_{Z_i},\end{equation} where for each $i$, $Z_i$ is an
algebraic set, $\one_{Z_i}$ is the characteristic function of
$X_i$, $f_i:Z_i\to X$ is a proper regular map, and $m_i$ is an
integer.

Clearly the sum of algebraically constructible functions is
algebraically constructible. The product of algebraically
constructible functions is algebraically constructible because
the fiber product of algebraic sets over $X$ is an algebraic set
over $X$: If $f_1:Z_1\to X$ and $f_2:Z_2\to X$ are proper
regular maps, then so is the fiber product $f:Z_1\times_X Z_2
\to X$, \begin{equation*} \begin{CD} Z_1\times_X Z_2 @>>> Z_2\\
@VVV @VV{f_2}V\\ Z_1 @>>f_1>      X \end{CD} \end{equation*}
where $ Z_1\times_X Z_2 = \{(z_1,z_2)\ |\ f_1(z_1) = f_2(z_2)\}$
and $f(z_1,z_2) = f_1(z_1) = f_2(z_2)$. Furthermore, for all
$x\in X$, $f\inv(x) = f_1\inv(x)\times f_2\inv(x)$. Therefore
$(f_{1*}\one_{Z_1})( f_{2*}\one_{Z_2}) = f_*\one_Z$.

Similarly, the pullback (\ref{pullback}) of an algebraically
constructible function by a regular map is algebraically
constructible. The pushforward (\ref{pushforward}) of an
algebraically constructible function by a proper regular map is
algebraically constructible, by the functoriality of pushforward
($(g\circ f)_* = g_*\circ f_*$).

The fact that the link (\ref{link}) of an algebraically
constructible function is algebraically constructible follows
from resolution of singularities and the fact that the link
operator commutes with pushforward ($\link f_* = f_*\link$).
Resolution of singularities implies that all the algebraic sets
$Z_i$ in (\ref{acf}) may be taken to be smooth and irreducible.
If $Z$ is a smooth algebraic set of dimension $d$, then
\begin{equation*} \link\one_Z=\begin{cases} 2\cdot\one_Z &
\text {$d$ odd}\\ 0 & \text{$d$ even}\end{cases}
\end{equation*} Thus if $\varphi = \sum m_i f_{i*}\one_{Z_i}$ as
in (\ref{acf}), \begin{align*} \link\varphi &= \link\sum m_i
f_{i*}\one_{Z_i}\\ &= \sum m_i f_{i*}\link\one_{Z_i}\\ &= 2\sum
m_i f_{i*}\one_{Z_i}, \end{align*} where the last sum is over
all $i$ such that $\dim Z_i$ is odd. This argument actually
gives the following stronger result \cite{mccpar2}.

\begin{thm}\label{halflink} If $\varphi$ is an algebraically
constructible function on the algebraic set $X$, then the values
of $\link\varphi$ are even, and $\half\link\varphi$ is an
algebraically constructible function. \end{thm}

This theorem is the key to defining Akbulut-King numbers using
constructible functions (section \ref{top} below). We call
$\hlink =\half\link$ the \emph{half-link} operator.

The constructible function $\varphi$ is \emph{Euler} if all the
values of $\link\varphi$ are even. The following examples show
that the sets of constructible functions, Euler constructible
functions, algebraically constructible functions, and strongly
algebraically constructible functions are all different. (More
examples can be found in \cite{mccpar2}.)

Let $X=\R^2$, and let $Q$ be the closed first quadrant,
\begin{equation*} Q = \{(x,y)\ |\ x\geq 0, y\geq 0\}.
\end{equation*} The constructible function $\one_Q$ is not
Euler, for $\link\one_Q = \one_O + \one_A + \one_B$, where $O$
is the origin, $A$ is the positive $x$-axis, and $B$ is the
positive $y$-axis.

The constructible function $2\cdot\one_Q$ is Euler, since it is
even, but it is not algebraically constructible. If it were
algebraically constructible, then by Theorem \ref{halflink} the
half-link $\hlink(2\cdot\one_Q) = \link\one_Q$ would be
algebraically constructible. But if $C$ is the $x$-axis, then
$(\one_C)(\link\one_Q) = \one_O + \one_A$, which is not Euler
and hence not algebraically constructible.

The constructible function $4\cdot\one_Q$ is not strongly
algebraically constructible, since the algebraic closure of the
first quadrant is the plane. But $4\cdot\one_Q$ is algebraically
constructible, \begin{equation*} 4\cdot\one_Q = f_*\one_{\R^2} +
g_*\one_{\R} + h_*\one_{\R} + \one_O,
\end{equation*} where $f(x,y) = (x^2,y^2)$, $g(x) = (x^2,0)$,
and $h(x) = (0,x^2)$.


\medskip \section{Real Algebra}\label{realalg}

The bridge between the topological and algebraic properties of
algebraically constructible functions  is given by the following
theorem \cite{parszaf1}, 
\cite{coskur2},  \cite{parszaf2}.

\begin{thm}\label{szafraniec1} Let $X$ be an algebraic subset of
$\R^n$ and let $f:Z\to X$ be a regular map.  Then there are real
polynomials $g_1,\ldots, g_s
\in \R[x_1,\ldots,x_n]$ such that for all $x\in X$,
\begin{equation}\label{presentation} \chi (f\inv (x)) =
\sum_{i=1}^s \sign g_i(x).
\end{equation}   \end{thm}

For a polynomial, or more generally for a regular function
$g$ on $X$, the projection $\pi:Z\to X$ from the  double cover
$Z = \{(x,t) \in X\times \R\ |\ t^2 = g(x)\}$ to $X$ satisfies
$\sign g = \pi_*\one_{Z} - \one_X$.  Thus the sign of
$g$ and hence any finite sum of signs of polynomials on
$X$ is algebraically constructible.  This gives the following
characterization of algebraically constructible functions. 

\begin{cor}\label{szafraniec2} Let $X$ be an algebraic set. 
Then $\varphi:X\to \Z $ is  algebraically constructible if and
only if $\varphi$ equals a finite sum of signs of polynomials on
$X$.
\end{cor}

Usually we suppose in the presentation
\eqref{presentation} that none of polynomials $g_i$ vanishes
identically on $X$.  Suppose, moreover, that $X$ is
irreducible.   Then the product $g = \prod g_i$ does not vanish
identically on $X$ and hence the zero set $W$ of 
$g$ is a proper algebraic subset of $X$.  Since $X$ is
irreducible, $\dim W < \dim X$.   By \eqref{presentation}, for
$x\in X\setminus W$, the Euler characteristic
$\chi (f\inv (x))$ is congruent mod 2 to $s$, the number of
polynomials $g_i$.  This gives  an alternative proof of Theorem
\ref{family}. We can go further and describe the Euler
characteristic of fibres of $f$ mod $4$: For $x\in X\setminus
W$,  
\begin{equation}\label{discriminant}
\chi (f\inv (x)) \equiv (s-1) + \sign g (x) \pmod 4.
\end{equation}  The existence of such $g$, called a
\emph{discriminant} of
$f$, was proved by Coste and Kurdyka in \cite{coskur2}. Although
the existence of a discriminant follows from
\eqref{presentation}, historically \eqref{discriminant} was a
prototype for Theorem \ref{szafraniec1}.   In general, a
polynomial or a rational function $g$ on $X$ satisfying
\eqref{discriminant} is not unique.  It is uniquely defined as
an element of the multiplicative group of non-zero rational
functions 
$\R(X)\setminus \{0\}$ on $X$ divided by the subgroup of those
functions which are generically of constant sign on
$X$, that is to say (by Hilbert's seventeenth problem), by
$\pm$ sums of squares of rational functions
(\emph{cf}.~\cite{coskur2}).  

Theorem \ref{szafraniec1} and Corollary \ref {szafraniec2}
provide an abstract link to the theory of quadratic forms.   Let
us consider a more concrete  example. Consider a finite regular
map $f:Z\to \R^n$,  where $Z\in \R^n\times \R$ is the zero set
of a polynomial without multiple factors,
\begin{equation*} P(x,z) = z^s + \sum_{i=0}^{s-1} a_i(x) z^i,
\end{equation*}
$x\in \R^n, z\in \R$,  and $f$ is induced by the projection on
the first factor. Since $f$ is finite, $\chi(f\inv (x))$ equals
the number of distinct real roots of the polynomial
$P_x(z)=P(x,z)$ of one real  variable.  By a classical theorem
of Hermite-Sylvester (see \emph{e.g.}~\cite{bcr} Proposition
6.2.6), the number of distinct real roots of $P_x$ equals the
signature of a symmetric matrix $Q$ of size $s$ with entries
polynomials in the $a_i$. The determinant of $Q$ equals the
discriminant
$\Delta (x)$ of $P_x$.  It is elementary to show that, if
$\Delta (x)$ is non-zero, then the number of distinct real roots
of $P_x$ is congruent to $s - 1 + \sign \Delta (x)$ mod $4$.  
This justifies the name ``discriminant''  for $g$ satisfying
\eqref{discriminant}.  If we want to compute the signature of
$Q$ we diagonalize it over the field of rational functions
$\R(x_1,\ldots, x_n)$. The signature of the matrix obtained, and
hence the signature of $Q$, equals the sum of signs of the
elements on the diagonal, which are rational functions.  The
sign of the rational function $p/q$ equals the sign of the
polynomial $pq$, in the complement of the zero set of the
denominator.  Thus we have shown the existence of polynomials
$g_i$ that satisfy \eqref{presentation} generically; that is,
for $x$ in the complement of a proper algebraic subset of
$\R^n$.    

In general, if $X$ is irreducible then we say that a function is
defined or a property holds \emph{generically} on $X$ if this is
so in the complement of a proper algebraic subset of $X$.  For
instance we call an integer-valued function on 
$X$ 
\emph{generically algebraically constructible} if it coincides 
with an algebraically constructible function  in the complement
of a proper algebraic subset of $X$.   By Corollary
\ref{szafraniec2}, 
 $\varphi$ is generically algebraically constructible on $X$ if
and only if it is  generically equal to the signature of a
quadratic form over the field of rational functions $\R(X)$. 
Thus for any regular map $f:Z\to X$ there is a quadratic form
$Q$ over $\R(X)$ such that
$\chi (f\inv (x))$  generically equals the signature of
$Q$.  There are, of course, many quadratic forms with this
property and we do not know whether, for arbitrary $f$, there is
a natural choice of such a  quadratic form.   The proofs of
Theorem \ref{szafraniec1} are all based on the  construction of
quadratic forms.  In
\cite{coskur2} the forms are constructed by means of Morse
Theory, in \cite{parszaf1} they are given by the Eisenbud-Levine
Theorem, and  in \cite{parszaf2} by a modern version of the
Hermite-Sylvester theorem.
 
Corollary \ref{szafraniec2} shows immediately that the sum and
the product of algebraically constructible functions are
algebraically constructible.  This corollary was used in  
\cite{parszaf1} to give another proof that
 the family of algebraically constructible functions is
preserved by the half-link operator (Theorem \ref{halflink}),
and by related topological operators such as  specialisation
(\emph{cf}.~\cite{mccpar2}, \cite{parszaf1}).  In general,
Corollary \ref{szafraniec2} allows us to use algebraic methods
to study algebraically constructible functions.

Recall that a \emph{basic} open semialgebraic subset of a real
algebraic set
$X\subset \R^n$ is a subset of the form
\begin{equation*}
\mathcal U(g_1,\ldots, g_s) = \{x\in X\ |\ g_1(x) >0 , \ldots,
g_s(x) >0 \},
\end{equation*}
 where $g_i$ are polynomials on $X$. For a polynomial $g$, the
function
$2\cdot\one_{\mathcal U(g)} = \sign g + \one_X -
\one_{g\inv (0)}$ is algebraically constructible, and hence so
is 
$2^s\one_{\mathcal U(g_1,\ldots, g_s)} = \prod 2\cdot
\one_{\mathcal U(g_i)}$.  By a theorem  of Br\"ocker and
Scheiderer  (\emph{cf}.~\cite{bcr} section 6.5), every basic
open semialgebraic subset of a real algebraic set of dimension
$d$ can be defined by at most
$d$ simultaneous strict polynomial inequalities; that is, we may
always choose $s\le d$.  This gives, as shown in \cite{mccpar2},
the following result. 

\begin{thm}\label{2d} Let $X$ be a real algebraic set of
dimension
  $d$.  Then any constructible function on $X$ divisible by
$2^d$ is
  algebraically constructible.  \end{thm}

For instance any constructible function on the plane $\R^2$
divisible by $4$ is algebraically constructible.  We do not know
an elementary proof of this fact.

As we have mentioned before, the first version of Theorem
\ref{costekurdyka} was proved by Coste \cite{cos1} using the
real spectrum.  
 This approach was later developed by Bonnard,  who introduced a
fan criterion for algebraically constructible functions
\cite{bon1}. Fans, which are subsets of the real spectrum, were
introduced by  Br\"ocker in order to study  quadratic forms.  

The fan criterion has proved to be a very powerful tool.  It
gives  (\emph{cf}.~\cite{bon1}) the following ``wall''
criterion.   Suppose that $X$ is nonsingular and compact,  
$\varphi:X\to \Z$, and $W$ is an algebraic subset of $X$ with
normal crossings such that $\varphi$ is constant on the
connected  components of $X\setminus W$.  For a smooth point
$w\in W$ we define $\partial_W\varphi (w)$ as the average of the
values of 
$\varphi$ on $X\setminus W$ computed on  both sides of $W$ at
$w$.  We may extend $\partial_W\varphi (w)$  arbitrarily to the
singular part of $W$.  Then  
$\varphi$ is generically algebraically constructible on $X$ if
and only if $\partial_W\varphi(w)$ is generically algebraically
constructible on $W$. This allows one to use induction on
dimension in working with generically algebraically
constructible functions.  

One may ask how many polynomials are necessary in order to
describe   a given algebraically constructible function as a sum
of  signs of polynomials.   In \cite{bon1} Bonnard gives two
bounds:    for a complete presentation and for a generic
presentation.  

\begin{thm}\label{bound}  Let $X$ be a real algebraic set of
dimension
$d$ and let $\varphi:X\to \Z$ be an algebraically constructible
function such that $\varphi (X) \subset [\delta -k, \delta +k]$,
$\delta \in\Z$, $k\in\N$.
\begin{itemize}
\item [(i)] Then  $\varphi$
  can be written as the sum of signs of at most $N'(d,k,\delta)$
  polynomials,  where $N'(d,k,\delta)$ equals
$2^{d-1}3k+|\delta|$ for $k$ even and 
  $2^{d-1}3(k-1)+2d+|\delta|$ for $k$ odd.
\item [(ii)] Suppose, moreover, that $X$ is irreducible.    
Then  $\varphi$ can be written generically  as the sum if signs
of at most $N(d,k,\delta)$ polynomials,  where $N(d,k,\delta)$
equals $2^{d-1}k+|\delta|$ for $k$ even and 
  $2^{d-1}(k-1)+1+|\delta|$ for $k$ odd.
\end{itemize}
\end{thm}

These bounds are sharp in the following sense.  Given $d, k,$ and
$\delta$ as above, there is an algebraic set $X$ of dimension $d$
and an  algebraically constructible function $\varphi'$,
resp.~$\varphi$,  on $X$ such that
$N'(d,k,\delta)$, resp.~$N(d,k,\delta)$,  is the minimal number
of polynomials necessary  to present $\varphi'$, resp.~to
present $\varphi$ generically,   as a sum of signs of
polynomials. Moreover, we may take $X=\R^d$.  Theorem
\ref{bound} is proved using 
 spaces of orderings.

For a fixed algebraically constructible function the bound given
by  Theorem \ref{bound} may not be sharp.  Let $X$ be
irreducible and  compact, and let $\varphi$ be an algebraically
constructible function on $X$
 constant on the connected  components of $X\setminus W$, where
$W\subset X$ is an algebraic subset with normal crossings.  In
\cite{bon4} Bonnard gives a recursive method for effectively
calculating the minimal number of polynomials representing
$\varphi$.

Using the real spectrum Pennaneac'h has defined \emph{algebraically
constructible chains} \cite{pen} and  algebraically
constructible  homology of real algebraic varieties.  An algebraically
constructible 
$k$-chain of a real algebraic variety $X$ is an oriented semialgebraic chain 
supported by a $k$-dimensional irreducible algebraic subset $V\subset X$,  
and ``weighted'' by a generically algebraically constructible function on $V$.  
The boundary operator is given by half of the standard boundary.  
Pennaneac'h has  proved that a constructible function is algebraically
constructible if and only if its characteristic Lagrangian cycle
is algebraically constructible.


\medskip \section{Topological Invariants}\label{top}

Using algebraically constructible functions and the link
operator we define local topological invariants which generalize
the Akbulut-King numbers. The vanishing of these invariants
gives necessary conditions for a topological space to be
homeomorphic to an algebraic set.

Algebraic sets are necessarily triangulable (\emph{cf.}
\cite{bcr}, 9.3.2), and by definition a triangulable space is
homeomorphic to a Euclidean simplicial complex, which is a
semialgebraic set. So without loss of generality we assume that
the spaces we consider are semialgebraic sets in Euclidean space.

If $X$ is a semialgebraic set, let $\cons(X)$ be the set of
constructible functions on $X$ (\ref{cf}). The set $\cons(X)$ is
a commutative ring with identity $\one_X$, and it is equipped
with a linear operator $\link$, the link operator, and a linear
integer-valued function $\varphi\mapsto
\int\varphi\dchi$, the Euler integral.

Of course $\cons(X)$ is not a topological invariant of $X$, but
the identity element, the link operator, and the Euler integral
are topological invariants in the following sense. Let $h:X'\to
X$ be a homeomorphism of semialgebraic sets. Then
$\one_X'=\one_{X}\circ h $. Let $\varphi\in \cons(Y)$ be such
that $\varphi'= \varphi\circ h \in \cons(X')$. Let $Y\subset X$
be a compact semialgebraic subset such that $Y'=h\inv(Y)$ is
also semialgebraic. Then \begin{align*}
\link(\varphi')&=(\link\varphi)\circ h,\\ \int_{Y'}\varphi'
\dchi&= \int_Y\varphi\dchi. \end{align*} For the elementary
proof see \cite{mccpar2}, appendix A.7. It follows that the
subring of $\cons(X)$ generated by $\one_X$ and $\link$ is a
topological invariant of $X$.

For a semialgebraic set $X$ let $\hlink(X)$ be the subring of
$\cons(X)\otimes\Q$ generated by $\one_X$ and the half-link
operator $\hlink = \half\link$. The ring $\hlink(X)$ is a
topological invariant of $X$. Theorem \ref{halflink} says that
if $X$ is an algebraic set then $\hlink(X)\subset \cons(X)$. In
other words, all of the functions obtained from $\one_X$ by the
arithmetic operations of sum, difference, and product, together
with the half-link operator, are \emph{integer-valued}.

So we have a method to produce topological obstructions for a
semialgebraic set $X$ to be homeomorphic to an algebraic set:
Find expressions $\varphi$ built from $\one_X$ using the
operations $+$, $-$, $\times$, and $\hlink$ such that $\varphi$
is not integer-valued. In \cite{mcc} (see also \cite{cos2}) this
method is explicitly illustrated for Akbulut and King's original
example of a compact 3-dimensional Euler space which is not
homeomorphic to an algebraic set. In \cite{mccpar2} the authors
show that this method reproduces the Akbulut-King invariants for
sets of dimension at most 3.

In \cite{mccpar3} we find all invariants produced by this method
for sets of dimension at most 4. The number of independent
invariants is enormous: $2^{29}-29$. But there is another
surprise. It follows from Corollary \ref{szafraniec2} that for
all algebraically constructible functions $\varphi$, the function
$\half(\varphi^4-\varphi^2)$ is algebraically constructible
(\cite{mccpar3}, Lemma 4.1). So a necessary condition for $X$ to
be algebraic is that all of the functions obtained from
$\one_X$ by the arithmetic operations of sum, difference, and
product, together with the half-link operator \emph{and} the
operator $\operatorname{P}(\varphi) =\half(\varphi^4-\varphi^2)$
are integer-valued. The total number of independent invariants
produced by our method taking all of these operators into
account is $2^{43}-43$.

The classification of these invariants is simplified by studying
them locally. The function $\varphi\in\cons(X)$ is Euler if and
only if, for every $x\in X$ and every $\e = \e(x)$ sufficiently
small, the restriction of $\varphi$ to $L_\e(x) = S(x,\e)\cap X$
has even Euler integral.  Now $\one_X|L _\e(x) =
\one_{L _\e(x)}$, and $\hlink(\varphi)|L _\e(x) = \varphi|L
_\e(x) - \hlink(\varphi|L _\e(x))$,  by \cite{mccpar3} 1.3(d).
Thus all functions obtained from $\one_X$ by the operations
$+$, $-$, $\times$, $\hlink$, and $\operatorname{P}$ are
integer-valued if and only if, for all $x\in X$, all functions
obtained from $\one_{L _\e(x)}$ by these operations are
integer-valued and have even Euler integral.

Akbulut and King obtain their numbers from cobordism invariants
of resolution towers. Here we give a cobordism-style description
of the Akbulut-King numbers which does not involve resolution
towers.

Let $n$ be a nonnegative integer, and let $\mathcal A_n$ be the
set of homeomorphism classes of compact real algebraic sets of
dimension at most $n$. If $X$ is a compact real algebraic set,
we say that $X$ is a \emph{boundary} if there exists a real
algebraic set $W$ and a point $w\in W$ such that $X$ is
homeomorphic to the link of $w$ in $W$. Let $\mathcal B_n$ be
the subset of $\mathcal A_n$ consisting of homeomorphism classes
of boundaries.

\begin{prop}\label{double} For all compact algebraic sets $X$,
the disjoint union $X\sqcup X$ is a boundary. \end{prop}
\begin{proof} Suppose $X\subset\R^n$ is given by the polynomial
equation $f(x_1,\dots,x_n)=0$ of degree $d$. Consider the
homogeneous polynomial $g(x_1,\dots,x_n,x_{n+1})$ of degree $d$
such that $g(x_1,\dots,x_n,1)=f(x_1,\dots,x_n)$. Let
$W\subset\R^{n+1}$ be given by $g(x_1,\dots,x_n,x_{n+1})=0$.
Then $X\sqcup X$ is homeomorphic to the link of the origin in
$W$.\end{proof}

Let $\mathcal A_n'$ be the free abelian group on the set
$\mathcal A_n$ modulo the subgroup generated by elements of the
form $[X]+[Y]-[X\sqcup Y]$, and let $\mathcal B_n'$ be the
subgroup of $\mathcal A_n'$ generated by $\mathcal B_n$. Let
$\mathcal V_n = \mathcal A_n'/\mathcal B_n'$. By Proposition
\ref{double}, $\mathcal V_n$ is a vector space over $\Z/2$.

The Akbulut-King numbers are additive under disjoint union, so
they define a linear map $a:\mathcal A_2'\to (\Z/2)^5$, $a([X])
= (\chi_2(X), a_0(X), a_1(X),a_2(X),a_3(X))$. By Theorem
\ref{akbulutking}, if $\dim X\leq 2$ then $X$ is a boundary if
and only if $a([X]) =0$. (If $a([X]) = 0$ then by Theorem
\ref{akbulutking} the cone on $X$ is homeomorphic to an
algebraic set.) Thus $a$ induces an injective linear map
\begin{equation*} \mathbf a:\mathcal V_2\to (\Z/2)^5.
\end{equation*} Akbulut and King show that $\mathbf a$ is an
isomorphism by constructing $2$-dimensional algebraic sets
$Y_0$,  $Y_1$, $Y_2$, and $Y_3$ such that $a(\R\mathbb P^2)$,
$a(Y_0)$, $a(Y_1)$, $a(Y_2)$, and $a(Y_3)$ are linearly
independent (\cite{akbking2}, p.~195--197). Thus $(\R\mathbb
P^2, Y_0,  Y_1, Y_2, Y_3)$ is a basis for $\mathcal V_2$.

On the other hand, the authors' constructible function
invariants \cite{mccpar2} define a linear map $b:\mathcal
A_2'\to (\Z/2)^5$, $b([X]) = (\chi_2(X), b_1(X),
b_2(X),b_3(X),b_4(X))$, where the $b_i$ are the mod 2 reductions
of the following Euler integrals. Let $\ahlink$ be the operator
defined by $\ahlink(\varphi) = \varphi -
\hlink(\varphi)$. Let $\alpha =\hlink\one_X$, $\beta =
\ahlink(\alpha^2)$, and $\gamma = \ahlink(\alpha^3)$. Then
$b_1(X) =\int\alpha\beta\dchi$, $b_2(X) =
\int\alpha\gamma\dchi$, $b_3(X) = \int\beta\gamma\dchi$, and
$b_4(X)=\int\alpha\beta\gamma\dchi$. The map $b$ induces a
linear isomorphism \begin{equation*} \mathbf b: \mathcal V_2\to
(\Z/2)^5. \end{equation*} The reader may verify that $\mathbf b$
is an isomorphism by computing $b(Y_i)$, $i=0, 1, 2, 3$. (The
computation of $b_1(Y_2)$ is given in \cite{mcc}.)

From this viewpoint the authors' results for $4$-dimensional
sets can be summarized as follows.  The invariants of
\cite{mccpar3} (Theorem 4.7) are additive under disjoint union,
 so they define a linear map $ c:\mathcal A_3'\to (\Z/2)^N$,
$N=2^{43}-43$. If $X$ is a boundary of dimension at most 3, then
$c([X])= 0$. Thus $c$ induces a linear map
\begin{equation*} \mathbf c:\mathcal V_3\to (\Z/2)^N.
\end{equation*} We show that $\mathbf c$ is surjective by
constructing a large set $\mathcal S$ of compact
$3$-dimensional Euler spaces with vanishing local Akbulut-King
numbers.  By Theorem \ref{akbulutking} every space in $\mathcal
S$ is homeomorphic to a real algebraic set. We make $\mathcal S$
so big that $\{ c([X])\ |\ X\in\mathcal S\}$ spans
$(\Z/2)^N$. Thus the dimension of $\mathcal V_3$ is at least
$2^{43}-43$. But we do not know whether $\mathbf c$ is an
isomorphism, because we do not know if $c([X])= 0$ implies that
$X$ is a boundary.

\begin{rem}\label{compactness} The one-point compactification of
a real algebraic set  is homeomorphic to a real algebraic set
(\emph{cf}.~\cite{akbking0},
\cite {bcr} Proposition 3.5.3).  Hence  a locally compact
topological space $X$ is homeomorphic to a real algebraic set if
and only if its one-point compactification
$\cX = X \cup \{\infty\}$ is homeomorphic to a real algebraic
set.  Suppose $X$ is a closed semialgebraic subset of $\R^n$
and  the local invariants constructed above vanish at every
$x\in X$.  Then they have to vanish at
$\infty \in \cX$.  Indeed,  suppose that $\hlink(X)\subset
\cons(X)$.  Of the operations that define
$\hlink(\cX)$, that is  $+$, $-$, $\times$, and $\hlink$, all but
$\hlink$ preserve $\cons(\cX)$.  Hence if, contrary to our
claim, $\hlink(\cX)\not\subset \cons(\cX)$ then there is
$\varphi \in
\hlink(\cX)$ such that $\hlink \varphi$ is not integer-valued. 
By our  assumption the restriction of $\hlink \varphi$ to $X$
has to be integer-valued, and we may write 
$\hlink \varphi = \psi + \frac 1 2 \one_{\infty}$ with
 $\psi \in \cons(\cX)$.  Consequently the Euler integral
$\int_{\cX} \hlink \varphi  \dchi \not\in \Z$. This contradicts
the identity 
\begin{equation*} \int_{\cX} \hlink \dchi = 0,
\end{equation*} which holds for any constructible function with
compact support   (\cite{mccpar2} Corollary 1.3).  The same
argument can be applied to the extended invariants since the
operator
$P(\varphi) = \half(\varphi^4-\varphi^2)$ preserves
 integer-valued functions.  

This remark can be also applied to Sullivan's invariant and
Akbulut and King's numbers.  This means that the one-point
compactification of a semialgebraic Euler set is Euler, and that
in Theorem \ref{akbulutking} it suffices to suppose that
$X$ is homeomorphic to a closed, not necessarily compact,
semialgebraic subset of $\R^n$ or, equivalently, that
$\cX$ is triangulable.
\end{rem}


\medskip \section{Generalizations and open questions}\label{ques}

The method of algebraically constructible functions may be
extended to more general classes of sets.  Fairly complete
results have  already been obtained for arc-symmetric
semialgebraic sets.  Much less is known for real analytic sets.

Arc-symmetric sets were introduced by Kurdyka in \cite{kur1} and
applied in \cite{kur2} to show that any injective regular
self-map of a real algebraic variety is surjective.   A
semialgebraic subset
$A$ of a real algebraic set $X$ is called \emph{arc-symmetric}
if for each real analytic arc
$\gamma : (-\varepsilon , \varepsilon) \to X$ such that
$\gamma ((-\varepsilon,0))\subset A$, the entire image
$\gamma ((-\varepsilon,\varepsilon))$ is contained in $A$.
Then, by the curve selection lemma, $A$ is closed in $X$. For
instance a connected component of a real algebraic subset of $X$
is arc-symmetric.  For more examples see \cite{kur1}.  

Let $X$ be a real algebraic set. A \emph{Nash constructible
function} on $X$ is an integer-valued function which can be
written as a finite sum
\begin{equation*} \varphi = \sum m_i
f_{i*}\one_{Z'_i},\end{equation*} where for each $i$, $Z'_i$ is
a connected component of an algebraic set $Z_i$, 
$f_i:Z_i\to X$ is a proper regular map, and $m_i$ is an
integer.  By \cite{mccpar3},  $A\subset X$ is arc-symmetric if
and only if $A$ is closed in $X$ and the characteristic function
of $A$  is Nash constructible.  The half-link of a Nash
constructible function is Nash constructible and hence for $A$
arc-symmetric $\hlink (A) \subset  \cons (A)$.   This shows that
every arc-symmetric set is Euler  and  Akbulut and King's
invariants vanish at each point of a $3$-dimensional
arc-symmetric set.  Thus, by Theorem
\ref{akbulutking} and Remark \ref{compactness}, every 
3-dimensional arc-symmetric semialgebraic set is homeomorphic to
a real algebraic set.  

Suppose, moreover, that $X$ is compact.  Then, as shown by
Bonnard in
\cite{bon3}, most of the algebraic results presented in section
\ref{realalg} above have their analogs for Nash constructible
functions.  For instance, a function on $X$ is Nash
constructible if and only if it is a finite sum of signs of
blow-Nash functions.  An analytic function on $X$  is called
\emph{Nash} if its graph is semialgebraic. It is called
\emph{blow-Nash} if it can be made Nash after composing with a
finite sequence of blowings-up. There is a fan criterion for
Nash constructible functions, and the minimum number of
blow-Nash functions required to represent a Nash constructible
function admits the same bounds as in Theorem \ref{bound}. 
These bounds are again sharp.   

If $X$ is compact then for every arc-symmetric semialgebraic
subset $A\subset X$ of dimension
$4$,  all $2^{43} -43$ local invariants constructed in section
\ref{top} vanish at each $x\in A$.  Without the assumption of
compactness of $X$ we do not know whether the operator
$P(\varphi)=  \half(\varphi^4-\varphi^2)$ preserves Nash
constructible functions on $X$, since we do not have a similar 
representation as finite sums of signs.  We may  conclude only
that the $2^{29}-29$ invariants constructed with the help of
$\hlink (A)$ vanish.

So far there has not been a similar study for real analytic
sets.  One natural possibility would be to define analytically
constructible functions on a real analytic set $X$ as  in
\eqref{acf} with the $f_i$ proper real analytic maps. The set of
such functions forms a ring.  By resolution of singularities
this ring is preserved by the half-link operator
$\hlink$.  Hence we conclude that  real analytic sets are Euler,
a result due to Sullivan \cite{sul}.  Akbulut and King's local
invariants vanish on real analytic sets, as well as the
$2^{29}-29$ invariants constructed in section \ref{top}. Again,
we do not know whether the operator
$P(\varphi)=  \half(\varphi^4-\varphi^2)$ preserves this family
of constructible functions.

\begin{question} \label{global} Does there exist a topological
proof that the operator
$P(\varphi)=  \half(\varphi^4-\varphi^2)$ preserves
algebraically constructible functions?  Does this follow from
the resolution of singularities, for instance by the method of
resolution towers? 
\end{question}

All of the constraints on the topology of real algebraic,
arc-symmetric, or real analytic sets obtained by our
constructions are local, and we do not know whether any
constraints of a global character exist.  

\begin{question} \label{noglobal} Let $X$ be a compact
triangulable topological space  and suppose that for each $x\in
X$ the link $\operatorname{lk}(x,X)$ is homeomorphic to a link
of a real algebraic set. Is $X$ homeomorphic to a real algebraic
set?
\end{question}

If $X$ is of dimension $\le 3$ then the answer is positive by
the results of Benedetti-Ded\`o and Akbulut-King (see section
\ref{ak}).   

\begin{question}\label{suff} Let $X$ be a compact triangulable
topological space
 of dimension $4$ such that for each $x\in X$ the  topological
invariants of 
 $\operatorname{lk}(x,X)$ constructed in section \ref{top}
vanish.   Is $X$ homeomorphic to
 a  real algebraic set?
\end{question}

In other words we ask whether the vanishing of our invariants is
sufficient for a set to be homeomorphic to a real algebraic
set.  

\begin{question}\label{links} Let $X$ be a compact real
algebraic set of dimension $3$ and suppose that all the
topological invariants of $X$ constructed in section \ref{top}
vanish, \emph{i.e.}~$c([X]) = 0$.  Is then $X$ homeomorphic to a
link of a real algebraic set?
\end{question}

Equivalently, we ask whether the linear map $ c:\mathcal A_3'\to
(\Z/2)^N$, introduced section \ref{top}, is injective.   A
positive answer to Question \ref{suff} answers positively
Question \ref{links}.  These two questions are equivalent if we
suppose a positive answer to Question \ref{noglobal} for sets of
dimension $4$.

One possible way to aproach these questions is to reinterpret
our invariants in terms of Akbulut and King's resolution towers
\cite{akbking2}. 

\begin{question}\label{restowers} Let $X$ be a compact
triangulable topological space
 of dimension $4$ such that  for each $x\in X$ the  topological
invariants
 of $\operatorname{lk}(x,X)$ constructed in section \ref{top}
vanish.   Does $X$ admit a topological  resolution tower in the
sense of Akbulut and King 
\cite{akbking2}? What extra properties does such a resolution
tower satisfy?
\end{question} 

Akbulut and King \cite{akbking2} show that a resolution tower can
be made algebraic if it satisfies some additional properties.
These properties are stronger that those given by the resolution
of singularities of real algebraic sets. Thus a positive answer
to Question \ref{restowers} would not necessarily guarantee a
solution to  Question \ref{suff}. Previte \cite{pre} has done
interesting work on obstructions for a 4-dimensional
space to have an algebraic resolution tower. 

Following the original idea of Sullivan one may study the links
of real algebraic sets using complexification. The vanishing of
Akbulut and King's numbers on  real algebraic sets can be proved
using a relation between complex conjugation and the monodromy
of the complex Milnor fibre of an ordered family of functions
\cite{mccpar1}.  

\begin{question} Can one construct our $2^{43}-43$ topological
invariants in dimension $4$ by means of complexification?  
\end{question}


\end{document}